\documentclass[12pt,a4paper,reqno]{amsart}

\usepackage{amsmath}
\usepackage{amssymb}
\usepackage[all]{xy}

\newtheorem{theorem}{Theorem}[section]
\newtheorem*{theorem*}{Theorem}
\newtheorem{lemma}[theorem]{Lemma}
\newtheorem{prop}[theorem]{Proposition}
\newtheorem{cor}[theorem]{Corollary}
\theoremstyle{remark}
\newtheorem{example}[theorem]{Example}
\newtheorem{remark}[theorem]{Remark}
\newtheorem{definition}[theorem]{Definition}

\newtheorem*{definition*}{Definition}

\def\N{{\mathbb N}}

\def\T{{\mathbb T}}

\def\C{{\mathbb C}}
\def\Z{{\mathbb Z}}

\newcommand{\isp}{\prec}

\newcommand{\lsp}{\operatorname{span}}

\newcommand{\Aut}{\operatorname{Aut}}
\numberwithin{equation}{section}
\newcommand{\Ind}{\operatorname{Ind}}
\newcommand{\id}{\operatorname{id}}
\begin{document}

\title{Corners of graph algebras} \date{\today}

\author{Tyrone Crisp}\address{School of Mathematical and Physical
Sciences, The University of Newcastle, Callaghan, NSW 2308,
Australia} \email{tyrone.crisp@studentmail.newcastle.edu.au}

\begin{abstract} It is known that given a directed graph $E$ and a subset $X$ of vertices, the sum $\Sigma_{v\in
X}P_v$ of vertex projections in the $C^*$-algebra of $E$ converges
strictly in the multiplier algebra to a projection $P_X$. Here we
give a construction which, in certain cases, produces a directed
graph $F$ such that $C^*(F)$ is isomorphic to the corner $P_X
C^*(E)P_X$. Corners of this type arise naturally as the fixed
point algebras of discrete coactions on graph algebras related to
labellings. We prove this fact, and show that our construction is
applicable to such a case whenever the labelling satisfies an
analogue of Kirchhoff's voltage law.
\end{abstract}


\maketitle

\section{Introduction}\label{introduction}

The $C^*$-algebra of a directed graph $E=(E^0,E^1,s,r)$ is a
universal object generated by Hilbert space operators satisfying
certain relations, where the relations reflect the path structure
of the graph \cite{BPRS,CK,FLR,KPR,KPRR,RS}. The vertices $v$ of
$E$ correspond to projections $P_v$ onto mutually orthogonal
subspaces, and the edges $e$ correspond to partial isometries
$S_e$ which map between these subspaces. Given any subset $X$ of
vertices, the sum $\Sigma_{v\in X}P_v$ converges strictly to a
projection, which we denote by $P_X$, in the multiplier algebra
$M(C^*(E))$ \cite{BPRS}. Corners of the type $P_X C^*(E)P_X$
associated to certain sets $X$ of vertices arise often in the
study of graph algebras: see \cite{BP,CG,D,DT,HSQLS,SzK}, and
Section \ref{coactions}.

It is in general very useful to be able to identify an abstract
$C^*$-algebra with a graph algebra. This is because the graphical
presentation encodes a great deal of structural information about
its associated $C^*$-algebra, and allows one to compute the
algebra's invariants via straightforward calculations
\cite{BHRS,DT2,HSPIS,RS}. So, given a directed graph $E$ and a
subset $X$ of vertices, it might be useful to realize the corner
$P_XC^*(E)P_X$ itself as the $C^*$-algebra of a directed graph. It
is to this goal that the first part of this paper is devoted. That
is, we shall give a construction which, in certain cases, produces
a graph for this corner from the graph $E$.

Probably the best known example of such a construction is the
procedure described in the literature as ``adding a tail to a
sink'', which is used to approximate the $C^*$-algebra of a graph
containing sinks as a full corner of the algebra of a graph
without sinks \cite[Lemma 1.2]{BPRS}. In a similar vein, if the
graph $E$ contains infinite-emitters, one may realize $C^*(E)$ as
a full corner of the $C^*$-algebra of a row-finite graph, via a
construction due to Drinen and Tomforde \cite{DT}. This
construction was generalized in \cite[Section 4]{BP}, and further
in \cite{CG}.

The conditions of \cite[Theorem 3.1]{CG} are in practice quite
limiting: for example, the theorem is not applicable if the
hereditary complement of the set $X$ contains loops, sinks or
infinite-emitters. In order to overcome these restrictions, the
approach taken here is substantially different to that of
\cite{BP,CG,DT}. The prototype for our present construction is
\cite[Section 2]{SzK}, where the graph $E$ was assumed to be
finite, and the set $X$ to consist of a single vertex. We
generalize and simplify this construction, and fix up a slight
error. Unfortunately, this new approach is not applicable to our
basic (nonunital) examples, adding a tail at a sink and the
desingularization of \cite{DT}. However, our construction applies
in particular to all unital graph algebras (Lemma \ref{extree}),
and it has been shown that any graph algebra can be approximated
as a direct limit of unital graph algebras \cite{RS}.

One context in which corners of graph algebras arise naturally is
as fixed point algebras of certain discrete coactions on graph
algebras. Indeed, the motivating example for this research was the
construction in \cite{HSQLS} of quantum lens spaces as the fixed
point algebras of certain actions of finite cyclic groups on
quantum spheres, by analogy with construction of the classical
lens spaces. The actions in question arise from labellings, in the
sense of Kumjian and Pask \cite{KP}, who showed that the crossed
product of a graph algebra $C^*(E)$ by such a group action is
itself isomorphic to the $C^*$-algebra of a directed graph, called
the skew product graph. The work of Kumjian and Pask was
generalized in \cite{DPR,KQR} to cover labellings of directed
graphs by discrete (not necessarily abelian) groups. Labellings of
this sort give rise to discrete group coactions, rather than the
compact group actions of \cite{KP}, but the realization of the
crossed product as the graph algebra of a skew product still
works. In Theorem \ref{fixedpoint} we show that, just as in the
case of the quantum lens spaces, the fixed point algebras of these
discrete coactions may be recovered as corners of the skew product
graph algebras, and then give a condition on labellings which
ensures that we may use the construction of Section \ref{corners}
to realize these corners as graph algebras.

\section{Preliminaries}\label{preliminaries}

We adopt the standard nomenclature of directed graphs and graph
algebras, as found in \cite{BPRS}, for example, with the following
additions:

\subsection*{Directed graphs} Let $E$ be a directed graph, and let $m,n\in\N\cup\{\infty\}$ be
such that $n\geq m$. If $\mu\in E^m$ and $\nu\in E^n$ are paths of
length $m$ and $n$ respectively, such that $\nu_i=\mu_i$ for all
$i=1,\ldots,m$, then we say that $\mu$ is an \emph{initial
subpath} of $\nu$, and write $\mu\isp\nu$.

Each finite path $\mu\in E^*$ gives a finite sequence
$s(\mu_1),r(\mu_1),r(\mu_2),\ldots,r(\mu)$ of vertices. The path
$\mu$ is called \emph{vertex-simple} if this sequence contains no
repeated vertices (i.e. if $\mu$ contains no loops). Similarly, an
infinite path $\nu\in E^\infty$ is called vertex-simple if its
corresponding right-, left- or bi-infinite sequence of vertices
contains no repetition. Each path of length zero (i.e. each
vertex) is also defined to be vertex-simple. A graph $E$ for which
$E^*\cup E^\infty$ contains only vertex-simple paths is called an
\emph{acyclic} graph. A graph $E$ for which $E^\infty$ contains no
vertex-simple paths is called a \emph{path-finite} graph.

If $E$ is a directed graph and $F$ a subgraph of $E$, then for
vertices $u,v\in E^0$ we write $u\geq_F v$ to mean that there is a
path $\mu\in F^*$ with $s(\mu)=u$ and $r(\mu)=v$. A subset
$X\subseteq E^0$ is said to be \emph{hereditary} if it has the
property that for all $v\in X$ and $u\in E^0$, $v\geq_E u$ implies
$u\in X$. For any subset $Y\subseteq E^0$ we shall denote by
$H_E(Y)$ the smallest hereditary subset of $E^0$ containing $Y$.
The set $H_E(Y)\setminus Y$ is referred to as the \emph{hereditary
complement} of $Y$ in $E$.

A subgraph $T$ of a directed graph $E$ is called a \emph{directed
subtree} of $E$ if it is acyclic and if $|T^1\cap r^{-1}(v)|\leq
1$ for each vertex $v\in T^0$ (that is, if each vertex in $T^0$
receives at most one edge in $T^1$). If $T$ is a directed subtree
of $E$, let $T^r$ denote the subset of $T^0$ consisting of those
vertices $v$ with $|T^1\cap r^{-1}(v)|=0$ (these vertices are
called the \emph{roots} of $T$). Let $T^l$ denote the subset of
$T^0$ consisting of those vertices $v$ with $|T^1\cap
s^{-1}(v)|=0$ (these vertices are called the \emph{leaves} of
$T$).

The concept of a directed subtree (in particular, a row- and
path-finite one) is central to our construction in Section
\ref{corners}, and the following lemma points out several basic
and useful facts about such graphs.

\begin{lemma}\label{T} Let $T$ be a row-finite, path-finite directed
subtree of a directed graph $E$. Then the following hold:
\begin{enumerate} \item For each $v\in T^0$ there exists a
unique path $\tau(v)$ in $T^*$ with source in $T^r$ and range $v$.
Then for $u,v\in T^0$, $v\geq_T u$ if and only if
$\tau(v)\isp\tau(u)$. \item For each $v\in T^0$ there exist at
most finitely many vertices $u\in T^0$ with $v\geq_T u$. \item For
each $v\in T^0$ there exists at least one $u\in T^l$ such that
$v\geq_T u$. \item Suppose $u,v\in T^0$ have $\tau(v)\isp\tau(u)$
and $u\neq v$. Then there exists a unique edge $e\in s^{-1}(v)\cap
T^1$ such that $\tau(v)e\isp\tau(u)$. If $f\in s^{-1}(v)\cap T^1$
satisfies $\tau(u)\isp\tau(v)f$ then $f=e$ and $\tau(v)e=\tau(u)$.
\end{enumerate}\end{lemma}

\begin{proof} (1) Fix $v\in T^0$. If $v\in T^r$ then $\tau(v)=v$.
If not, then $v$ receives exactly one edge $e_1\in T^1$. If
$s(e_1)$ is in $T^r$ then $\tau(v)=e_1$. If not, then $s(e_1)$
receives exactly one edge $e_2\in T^1$. As $T$ is path-finite we
get a path $\tau(v)=e_n e_{n-1}\ldots e_1$ with source in $T^r$
after finitely many iterations of this construction. Uniqueness of
$\tau(v)$ follows from the fact that each $u\in T^0$ receives at
most one edge, and this same fact gives the equivalence $v\geq_T
u\iff \tau(v)\isp\tau(u)$.

(2) Suppose $v\in T^0$ is such that infinitely many such $u$
exist. As $T$ is row-finite, there is an edge $e_1\in
s^{-1}(v)\cap T^1$ such that $e_1$ is the first edge in infinitely
many distinct paths in $T^*$ (by the pigeonhole principle). We may
apply this same argument to the vertex $r(e_1)$, giving an edge
$e_2\in s^{-1}(r(e_1))\cap T^1$ such that $e_1 e_2$ is an initial
subpath of infinitely many distinct paths in $T^*$. Continuing
this construction gives a path $e_1 e_2\ldots\in T^\infty$, which
must be vertex-simple because $T$ is acyclic. This contradicts the
assumption that $T$ is path-finite, proving the claim.

(3) Fix $v\in T^0$. If $v\in T^l$ then we are done. Otherwise
choose $e_1\in s^{-1}(v)\cap T^1$. If $r(e)\in T^l$ then we are
done; otherwise find $e_2\in s^{-1}(r(e_1))\cap T^1$. This
construction must terminate after finitely many iterations,
because $T$ is path-finite.

(4) Fix $u,v\in T^0$ with $\tau(v)\isp\tau(u)$ and $u\neq v$.
Clearly there exists an edge $e\in s^{-1}(v)\cap T^1$ such that
$\tau(v)e\isp\tau(u)$. Suppose $e'$ is another such edge. Then
$r(e)=r(e')\in T^0$, contradicting that each edge in $T^0$
receives at most one edge. Hence $e$ is unique. Now suppose $f\in
s^{-1}(v)\cap T^1$ has $\tau(u)\isp\tau(v)f$. Since
$\tau(v)\isp\tau(u)$ we must then have $\tau(u)=\tau(v)f$, so
$\tau(v)f\isp\tau(u)$ and $f=e$ by uniqueness of $e$.
\end{proof}

\subsection*{Graph $C^*$-algebras} A Cuntz-Krieger $E$-family is a set
$\{P_v,S_e:v\in E^0,e\in E^1\}$ of operators on a Hilbert space
such that the elements $P_v$ are mutually orthogonal projections
and the elements $S_e$ are partial isometries with mutually
orthogonal ranges, satisfying the following relations:

\begin{enumerate}\item[(CK1)] $S_e^* S_e=P_{r(e)}$ for all $e\in E^1$;
\item[(CK2)] $S_e S_e^*\leq P_{s(e)}$ for each $e\in E^1$;
\item[(CK3)] $P_v=\sum_{e\in s^{-1}(v)} S_e S_e^*$ for each $v\in
E^0$ with $0<|s^{-1}(v)|<\infty$.\end{enumerate}

The $C^*$-algebra of $E$, denoted $C^*(E)$, is defined to be the
universal $C^*$-algebra generated by a Cuntz-Krieger $E$-family.

For any subset $X\subseteq E^0$, the sum $\sum_{v\in X} P_v$
converges strictly to a projection $P_X$ in $M(C^*(E))$
\cite[Lemma 1.1]{BPRS}.

\section{Corners of directed graphs} \label{corners}

In this section we describe our procedure for constructing a graph
for the corner $P_XC^*(E)P_X$ of a graph algebra $C^*(E)$
associated to a vertex set $X$. This construction is given in the
following definition, and its relation to $P_XC^*(E)P_X$ is shown
in Theorem \ref{main}.

\begin{definition}\label{corner} Let $E$ be a directed graph, $X\subseteq E^0$, and let
$T$ be a row-finite, path-finite directed subtree of $E$ with
$T^r=X$ and $T^0=H_E(X)$. Define a new directed graph, denoted
$E(T)$ and called the \emph{$T$-corner of $E$}, as follows:
\begin{align*}& E(T)^0:=T^0\setminus\{v\in T^0:\emptyset\neq
s^{-1}(v)\subseteq T^1\}\\ &E(T)^1:=\{e_u:e\in
s^{-1}(T^0)\setminus T^1,u\in E(T)^0,r(e)\geq_T u\}\\
&s(e_u)=s(e),\ r(e_u)=u.\end{align*}
\end{definition}

\begin{example} Suppose $X$ is a hereditary subset of $E^0$. We then have $T^0=X=T^r$, and
since no root may receive an edge in $T$ we infer that $T^1$ is
empty. Thus each $v\in T^0$ is either a sink, or emits an edge
which does not belong to $T^1$; this implies that $E(T)^0=T^0=X$.
Furthermore, for each edge $e$ with source in $T^0$, and each
vertex $u\in T^0$, $r(e)\geq_T u$ if and only if $r(e)=u$, because
$T^*=T^0$. Hence $E(T)^1=\{e_{r(e)}:s(e)\in X\}$, where each
$e_{r(e)}$ has the same range and source as $e$. Thus $E(T)$ is
nothing but the graph $(X,s^{-1}(X),s,r)$.
\end{example}

\begin{example}\label{O_2} Let $E$
be the graph \[\xygraph{{v_0}="v":_*{e_1}
@/_1pc/[dl(2)]{v_1}="dl":_*{e_2}@/_1pc/[r(4)]{v_2}="dr":_*{e_0}@/_1pc/"v":_-*{f_2}
@/_/"dr":_-*{f_1}@/_/"dl":_-*{f_0}@/_/"v"}\]

\noindent and let $T^0=E^0$, $T^1=\{e_1,f_2\}$. $T$ is a row- and
path-finite directed subtree of $E$ with root set $X=\{v_0\}$,
such that $T^0=H_E(X)$. The vertex $v_0$ is not a sink, and each
edge with source $v_0$ belongs to $T^1$, so $v_0\not\in E(T)^0$.
On the other hand, both $v_1$ and $v_2$ emit edges in $E$ which
are not part of $T$ (for example, $f_0$ and $e_0$ respectively),
so both belong to $E(T)^0$.

Now constructing the edge set $E(T)^1$, we consider in turn each
edge in $s^{-1}(T^0)\setminus T^1=\{e_0,e_2,f_0,f_1\}$; let us
start with $e_0$. The range $v_0$ of $e_0$ satisfies $v_0\geq_T
v_1$, because $e_1$ is a path in $T^*$ with source $v_0$ and range
$v_1$. Since $v_1$ belongs to $E(T)^0$ there will be an edge
$(e_0)_{v_1}$ in $E(T)^1$, with $s((e_0)_{v_1})=s(e_0)=v_2$ and
$r((e_0)_{v_1})=v_1$. Similarly, there will be an edge
$(e_0)_{v_2}$ with source $v_2$ and range $v_2$. Notice that,
although $v_0=r(e_0)\geq_T v_0$, there is no edge $(e_0)_{v_0}$
because $v_0\not\in E(T)^0$. Considering the remaining edges
$e_2,f_0,f_1\in s^{-1}(T^0)\setminus T^1$ in a similar way, we
obtain that
$E(T)^1=\{(e_0)_{v_1},(e_0)_{v_2},(e_2)_{v_2},(f_0)_{v_1},(f_0)_{v_2},(f_1)_{v_1}\}$,
so $E(T)$ is the following graph: \[\xygraph{
{v_1}="v_1":_{(f_0)_{v_1}}@(ul,dl)"v_1":^{(f_0)_{v_2}}@/^3pc/[r(2)]{v_2}="v_2"
:_{(e_0)_{v_2}}@(dr,ur)"v_2":^{(e_0)_{v_1}}@/^3pc/"v_1":^{(e_2)_{v_2}}@/^1pc/"v_2":^{(f_1)_{v_1}}@/^1pc/"v_1"}
\]
\end{example}

There will often be more than one choice of subtree with the
desired properties, giving nonisomorphic graphs $E(T)$:

\begin{example}\label{pqr} Let $p,q,r$ be positive integers, and let $E$ be as shown:
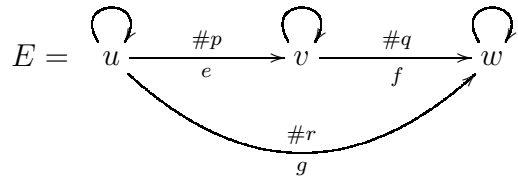
\begin{equation*}E=\xygraph{{u}:@(ul,ur)"u":^{\#p}_e[r(2)]{v}:
@(ul,ur)"v":^{\#q}_f[r(2)]{w}:@(ul,ur)"w",
"u":@/_3pc/^{\#r}_g"w"}\end{equation*}

\noindent Here the label ``$\#n$'' above an arrow indicates that
that arrow represents $n$ edges, and a label of ``$x$'' below an
arrow means that we will distinguish one of those edges and call
it $x$. Then the subgraph $T_1$ with $T_1^0=E^0$ and
$T_1^1=\{e,g\}$ is a finite directed subtree of $E$ with root
$X_1=\{u\}$, satisfying $T_1^0=H_E(X_1)$. For this $T_1$, the
construction of Definition \ref{corner} gives $E(T_1)\cong E$. On
the other hand, let $T_2$ be the finite subtree $T_2^0=E^0$,
$T_2^1=\{e,f\}$. This tree also has root set $X_2=X_1=\{u\}$ and
$T_2^0=H_E(X_2)$, but now the construction gives the following
graph:
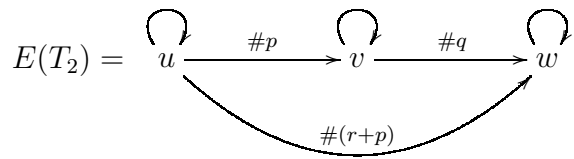
\begin{equation*}E(T_2)=\xygraph{{u}:@(ul,ur)"u":^{\#p}[r(2)]{v}:
@(ul,ur)"v":^{\#q}[r(2)]{w}:@(ul,ur)"w",
"u":@/_3pc/^{\#(r+p)}"w"}\end{equation*}
\end{example}

\begin{theorem}\label{main} Let $E$, $X$, $T$ and $ E(T)$ be as
in Definition \ref{corner}. Then $C^*(E(T))\cong
P_{X}C^*(E)P_{X}$.
\end{theorem}

Notice that the important properties of the tree are row- and
path-finiteness, as one can find a subtree with the other
property, for any root set  $X$, using an inductive construction
as in the proof of Lemma \ref{extree}. As an aside, the following
lemma indicates the scope of Theorem \ref{main}:

\begin{lemma}\label{extree} If $H_E(X)\setminus X$ is finite, then there is a subtree
with the desired properties. If $X$ is finite and $H_E(X)$
infinite, then there is no such subtree. \end{lemma}

\begin{proof} First suppose $H_E(X)\setminus X$ is finite. For
each $v\in H_E(X)$ let $d(v)$ be the length of a shortest path in
$E^*$ with source in $X$ and range $v$ (such a path exists because
$H_E(X)$ is hereditary). Let $T^0=H_E(X)$ and construct the edge
set $T^1$ recursively as follows. For each $n\in\N$ and each $v\in
H_E(X)\setminus X$ with $d(v)=n$ choose one edge $e_v\in E^1$ such
that $r(e_v)=v$ and $d(s(e_v))=n-1$. Let $T^1=\{e_v:v\in
H_E(X)\setminus X\}$. Then the subgraph $T$ of $E$ is row- and
path-finite by finiteness of $H_E(X)\setminus X$. On the other
hand, suppose $X$ is finite and $H_E(X)$ infinite, and suppose $T$
is a directed subtree of $E$ with roots $X$ and vertex set
$H_E(X)$. By Lemma \ref{T}(1) and the pigeonhole principle, there
must be a vertex $v\in X$ such that $v\geq_T u$ for infinitely
many vertices $u\in H_E(X)$. Hence, by part (2) of Lemma \ref{T},
$T$ cannot be row- and path-finite.
\end{proof}

The proof of Theorem \ref{main} will proceed in three main steps:
first we find a Cuntz-Krieger family for $E(T)$ inside $C^*(E)$,
so that the universal property of $C^*(E(T))$ gives a homomorphism
$\phi:C^*(E(T))\to C^*(E)$. Next we show that this $\phi$ is
injective, using the gauge-invariant uniqueness theorem
\cite[Corollary 1.4]{SzCKU}. Finally we show that the range of
$\phi$ is equal to $P_{X}C^*(E)P_{X}$ using an inductive argument.

For the first step, let $\{P_v,S_e\}$ be the canonical
Cuntz-Krieger generators of $C^*(E)$. For each $v\in T^0$ let
$\tau(v)\in T^*$ be the path given by part (1) of Lemma \ref{T}
(in particular, for $v\in X$, $\tau(v)=v$). Now for each $v\in
T^0$, define
\[Q_v:=S_{\tau(v)}S_{\tau(v)}^*-\sum_{e\in T^1\cap
s^{-1}(v)}S_{\tau(v)e}S_{\tau(v)e}^*.\] Since $T$ is row-finite,
this sum is finite and each $Q_v$ is an element of $C^*(E)$. The
relations (CK1)--(CK3) in $C^*(E)$ imply that each $Q_v$ is a
projection. These projections will correspond to the vertex
projections of $E(T)$, and we shall need to know that they are
nonzero:

\begin{lemma}\label{Q_v} For each $v\in T^0$, $Q_v=0$ if and only if $\emptyset\neq s^{-1}(v)\subseteq
T^1$. Also,
\begin{equation}\label{Q_v1}S_{\tau(v)}S_{\tau(v)}^*=\sum_{u\in
T^0,\,v\geq_T u}Q_u.\end{equation}\end{lemma}

\begin{proof} For the first claim, first suppose $\emptyset\neq s^{-1}(v)\subseteq
T^1$. The subgraph $T$ is row-finite, and each edge with source
$v$ belongs to $T^1$, so $0<|s^{-1}(v)|<\infty$. The Cuntz-Krieger
relation (CK3) in $C^*(E)$ then gives $P_v=\sum_{e\in
s^{-1}(v)}S_e S_e^*$, so we have \begin{align*}
Q_v&=S_{\tau(v)}S_{\tau(v)}^*-\sum_{e\in
s^{-1}(v)}S_{\tau(v)e}S_{\tau(v)e}^*\qquad\text{(since $T^1\cap s^{-1}(v)=s^{-1}(v)$)}\\
&=S_{\tau(v)}P_vS_{\tau(v)}^*-S_{\tau(v)}\left(\sum_{e\in
s^{-1}(v)}S_e S_e^*\right)S_{\tau(v)}^*=0.\end{align*} Conversely,
if $v$ is a sink in $E$ then $Q_v=S_{\tau(v)}S_{\tau(v)}^*\neq 0$.
If $v$ emits an edge $f\in E^1\setminus T^1$ then the relations
(CK1)--(CK3) in $C^*(E)$ imply that $S_{\tau(v)f}S_{\tau(v)f}^*$
is a subprojection of $S_{\tau(v)}S_{\tau(v)}^*$ orthogonal to
$\sum_{e\in T^1\cap s^{-1}(v)}S_{\tau(v)e}S_{\tau(v)e}^*$, so
$Q_v\geq S_{\tau(v)f}S_{\tau(v)f}^*\neq 0$.

For the second claim, first notice that the sum is finite by part
(2) of Lemma \ref{T}. For each vertex $v\in T^0$, let $c(v)$ be
the number of elements in the set $\{u\in T^0:v\geq_T u\}$. The
formula \eqref{Q_v1} will be derived by induction on $c(v)$. For
the basis step, note that if $c(v)=1$ then $T^1\cap
s^{-1}(v)=\emptyset$, so $Q_v=S_{\tau(u)}S_{\tau(u)}^*$ as
desired. For $n\in\N$, suppose the formula \eqref{Q_v1} holds for
all $w\in T^0$ with $c(w)\leq n-1$, and let $v\in T^0$ have
$c(v)=n$. Now
\begin{equation}\label{Q_v2}S_{\tau(v)}S_{\tau(v)}^*=Q_v+\sum_{e\in T^1\cap
s^{-1}(v)}S_{\tau(v)e}S_{\tau(v)e}^*,\end{equation} and for each
$e\in T^1\cap s^{-1}(v)$ we have $\tau(v)e=\tau(r(e))$ and
$c(r(e))<c(v)$, so $S_{\tau(v)e}S_{\tau(v)e}^*=\sum_{u\in
T^0,\,r(e)\geq_T u}Q_u$ by the inductive hypothesis. Substituting
this into \eqref{Q_v2} gives the formula \eqref{Q_v1} for the
vertex $v$.
\end{proof}

Now for $e_u\in E(T)^1$ define
$T_{e_u}:=S_{\tau(s(e))e}S_{\tau(r(e))}^*Q_u$.

\begin{prop}\label{CK} The family
$\{Q_v,T_{e_u}:v\in E(T)^0,\,e_u\in E(T)^1\}$ is a Cuntz-Krieger
family for the graph $ E(T)$.
\end{prop}

\begin{proof} The proof of this proposition requires some technical manipulations
of the relations (CK1)--(CK3), but is theoretically
straightforward. Lemma \ref{Q_v} implies that for each $v\in
E(T)^0$, $Q_v$ is a nonzero projection. To see that they are
mutually orthogonal, first notice that for each $v$, $Q_v$ is a
subprojection of $S_{\tau(v)}S_{\tau(v)}^*$. Suppose $v$ and $w$
are distinct elements of $ E(T)^0$ such that $Q_vQ_w\neq 0$. We
must have $S_{\tau(v)}^*S_{\tau(w)}\neq 0$, and hence one of
$\tau(v)$ and $\tau(w)$ is an initial subpath of the other (this
implication is a consequence of the fact that the $S_e$ have
mutually orthogonal ranges). Assume, without loss of generality,
that $\tau(w)\prec\tau(v)$, and let $f\in T^1\cap s^{-1}(w)$ be
the edge given by Lemma \ref{T}(4). Then
\[\sum_{e\in T^1\cap
s^{-1}(w)}S_{\tau(v)}^*S_{\tau(w)e}S_{\tau(w)e}^*=S_{\tau(v)}^*S_{\tau(w)f}S_{\tau(w)f}^*,\]
because $f$ is the unique edge in $T^1\cap s^{-1}(w)$ with the
property that $\tau(w)f\isp\tau(v)$. Now
$S_{\tau(v)}^*S_{\tau(w)f}S_{\tau(w)f}^*=S_{\tau(v)}^*$, and thus
\begin{align*}Q_vQ_w&=Q_vS_{\tau(v)}S_{\tau(v)}^*\left(S_{\tau(w)}S_{\tau(w)}^*-\sum_{e\in T^1\cap
s^{-1}(w)}S_{\tau(w)e}S_{\tau(w)e}^*\right)\\
&=Q_v\left(S_{\tau(v)}S_{\tau(v)}^*-S_{\tau(v)}S_{\tau(v)}^*\right)\\&=0.\end{align*}
Hence $Q_vQ_w\neq 0$ if and only if $v=w$.

Turning our attention to the $T_{e_u}$, fix $e_u\in E(T)^1$. By
definition of $ E(T)^1$ we must have $\tau(r(e))\prec\tau(u)$, so
$Q_u\leq S_{\tau(u)}S_{\tau(u)}^*\leq
S_{\tau(r(e))}S_{\tau(r(e))}^*$. Therefore
\begin{align*}
T_{e_u}^*T_{e_u}&=Q_uS_{\tau(r(e))}(S_{\tau(s(e))e}^*S_{\tau(s(e))e})S_{\tau(r(e))}^*Q_u
=Q_uS_{\tau(r(e))}P_{r(e)}S_{\tau(r(e))}^*Q_u\\&=Q_u(S_{\tau(r(e))}S_{\tau(r(e))}^*)Q_u
=Q_u.\end{align*} Thus the $T_{e_u}$ are nonzero partial
isometries with $T_{e_u}^*T_{e_u}=Q_{r(e_u)}$. To see that they
have mutually orthogonal ranges, take $e_u$ and $f_v$ in $ E(T)^1$
and suppose $T_{e_u}^*T_{f_v}\neq 0$. Now
\begin{equation}\label{mutorth}T_{e_u}^*T_{f_v}=Q_u
S_{\tau(r(e))}S_{\tau(s(e))e}^*S_{\tau(s(f))f}S_{\tau(r(f))}^*Q_v,\end{equation}
and in order for this product to be nonzero we must have either
$\tau(s(f))f\prec\tau(s(e))e$ or $\tau(s(e))e\prec\tau(s(f))f$.
Since neither $e$ nor $f$ belongs to $T^1$ (so that neither may be
part of any $\tau(w)$), this implies that
$\tau(s(e))e=\tau(s(f))f$, and so $e=f$. Putting $e=f$ in
\eqref{mutorth} gives
\[T_{e_u}^*T_{f_v}=Q_uS_{\tau(r(e))}S_{\tau(r(e))}^*Q_v=Q_uQ_v,\]
and in order for this product to be nonzero we must have $u=v$.
Thus $e_u=f_v$.

For the inequality $T_{e_u}T_{e_u}^*\leq Q_{s(e_u)}$, we calculate
\begin{align*}S_{\tau(s(e))e}^*Q_{s(e)}&=S_{\tau(s(e))e}^*\left(S_{\tau(s(e))}S_{\tau(s(e))}^*-\sum_{f\in
T^1\cap s^{-1}(s(e))}S_{\tau(s(e))f}S_{\tau(s(e))f}^*\right)\\
&=S_{\tau(s(e))e}^*S_{\tau(s(e))e}S_{\tau(s(e))e}^*-0=S_{\tau(s(e))e}^*,\end{align*}
since $e\not\in T^1$ implies that $\tau(s(e))e$ is not an initial
subpath of any $\tau(s(e))f$ for $f\in T^1$. Thus
\[T_{e_u}T_{e_u}^*Q_{s(e_u)}=T_{e_u}Q_vS_{\tau(r(e))}(S_{\tau(s(e))e}^*Q_{s(e)})
=T_{e_u}Q_vS_{\tau(r(e))}S_{\tau(s(e))e}^*=T_{e_u}T_{e_u}^*.\]

To prove the remaining identity (CK3), we need to know the
following fact about singular vertices in $E(T)$:

\begin{lemma}\label{sing} Each edge $e$ in $E^1\setminus T^1$ with
$s(e)\in T^0$ gives at least one edge in $ E(T)$ with source
$s(e)$. In particular, if $v\in T^0$ is a singular vertex of $E$
then $v$ is a singular vertex of $ E(T)$.\end{lemma}

\begin{proof} Let $e$ be an edge in $E^1\setminus T^1$ with
$s(e)\in T^0$, and let $s(e)=v$, $r(e)=u$. Since $T^0$ is a
hereditary subset of $E^0$ we must have $u\in T^0$. By part (3) of
Lemma \ref{T}, there exists at least one vertex $u'\in T^0$ with
$u\geq_T u'$ and $s^{-1}(u')\cap T^1=\emptyset$. Then by
definition of $ E(T)$ we have $u'\in E(T)^0$, and there is an edge
$e_{u'}$ in $ E(T)^1$ with source $v$.

For the second claim, suppose $v\in T^0$ is a sink in $E$. Then
$\emptyset=s^{-1}(v)$ gives $v\in E(T)$, and since there is no
edge in $E^1$ with source $v$ there is no edge $e_u$ in $ E(T)^1$
with source $v$. Hence $v$ is a sink in $E(T)$. On the other hand,
suppose $v\in T^0$ emits infinitely many edges. Since $T$ is
row-finite, infinitely many of these edges must belong to
$E^1\setminus T^1$. Each of these edges gives at least one edge in
$ E(T)$ with source $v$ by the preceding paragraph, so $v$ is an
infinite-emitter in $E(T)$.\end{proof}

Now suppose $v\in E(T)^0$ is nonsingular in $ E(T)$. Then $v$ is
nonsingular in $E$ by the preceding lemma, so the Cuntz-Krieger
relation (CK3) in $C^*(E)$ gives $P_v=\sum_{e\in
s^{-1}(v)}S_eS_e^*$. Now
\begin{align} Q_v&=S_{\tau(v)}S_{\tau(v)}^*-\sum_{e\in T^1\cap
s^{-1}(v)}S_{\tau(v)e}S_{\tau(v)e}^*=S_{\tau(v)}P_vS_{\tau(v)}^*-S_{\tau(v)}\left(\sum_{e\in
T^1\cap s^{-1}(v)}S_eS_e^*\right)S_{\tau(v)}^*\notag\\
&=S_{\tau(v)}\left(P_v-\sum_{e\in T^1\cap
s^{-1}(v)}S_eS_e^*\right)S_{\tau(v)}^*=\sum_{e\in
s^{-1}(v)\setminus
T^1}S_{\tau(v)e}S_{\tau(v)e}^*.\label{Q_vsum}\end{align}

Fix an edge $e\in s^{-1}(v)\setminus T^1$. This edge gives one
edge $e_u$ in $E(T)$ with source $v$ for each vertex $u\in E(T)^0$
with $r(e)\geq_T u$. The formula \eqref{Q_v1} of Lemma \ref{Q_v}
gives \begin{align*}
S_{\tau(v)e}S_{\tau(v)e}^*&=S_{\tau(v)e}P_{r(e)}^3S_{\tau(v)e}^*=S_{\tau(v)e}S_{\tau(r(e))}^*(S_{\tau(r(e))}S_{\tau(r(e))}^*)^2S_{\tau(r(e))}S_{\tau(v)e}^*\\
&=S_{\tau(v)e}S_{\tau(r(e))}^*(S_{\tau(r(e))}S_{\tau(r(e))}^*)\left(S_{\tau(v)e}S_{\tau(r(e))}^*(S_{\tau(r(e))}S_{\tau(r(e))}^*)\right)^*\\
&=\left(S_{\tau(v)e}S_{\tau(r(e))}^*\left(\sum_{u\in
T^0,\,r(e)\geq_T
u}Q_u\right)\right)\left(S_{\tau(v)e}S_{\tau(r(e))}^*\left(\sum_{u\in
T^0,\,r(e)\geq_T u}Q_u\right)\right)^*\\
&=\left(\sum_{u\in  E(T)^0,\,r(e)\geq_T
u}T_{e_u}\right)\left(\sum_{u\in  E(T)^0,\,r(e)\geq_T u}
T_{e_u}^*\right).\end{align*} Since for $u\neq u'$ we have
$T_{e_u}T_{e_{u'}}^*=0$, this product expands as
\[S_{\tau(v)e}S_{\tau(v)e}^*=\sum_{u\in  E(T)^0,\,r(e)\geq_T
u}T_{e_u}T_{e_u}^*.\] Substituting into \eqref{Q_vsum} now gives
the Cuntz-Krieger identity $Q_v=\sum_{s(e_u)=v}T_{e_u}T_{e_u}^*$,
and this final identity completes the proof of the proposition.
\end{proof}

Now the universal property of $C^*(E(T))$ gives a $*$-homomorphism
$\phi:C^*( E(T))\to C^*(E)$ which maps each canonical generator of
$C^*( E(T))$ to its corresponding element of the family
$\{Q_v,T_{e_u}\}$. The following two propositions show that $\phi$
is injective and has range $P_{X}C^*(E)P_{X}$.

\begin{prop}\label{inj} The map $\phi$ defined above is
injective.\end{prop}

\begin{proof} Arguing as in \cite[Section 1]{BPRS}, the universal property of $C^*(E)$ implies that there exists
an action $\sigma:\T\to\Aut(C^*(E))$ given on generators by
$\sigma_t(P_v)=P_v$ for all $v\in E^0$,
and\begin{equation*}\sigma_t(S_e)=\begin{cases}S_e&\text{for }e\in
T^1\\tS_e &\text{otherwise.}\end{cases}\end{equation*} This action
does not move any $S_{\tau(v)}$ for $v\in E(T)^0$, and hence does
not move any $Q_v$ either. For $e_u\in E(T)^1$ we have $e\in
E^1\setminus T^1$, and so for $t\in\T$,
\[\sigma_t(T_{e_u})=\sigma_t(S_{\tau(s(e))}S_eS_{\tau(r(e))}^*Q_v)=S_{\tau(s(e))}tS_eS_{\tau(r(e))}^*Q_v=tT_{e_u}.\]
Thus if $\gamma$ denotes the gauge action on $C^*(E(T))$ we have
$\phi\circ\gamma=\sigma\circ\phi$, and all $Q_v$ are nonzero, so
the gauge-invariant uniqueness theorem \cite[Corollary 1.4]{SzCKU}
implies that $\phi$ is injective.
\end{proof}

\begin{prop}\label{surj} $\phi(C^*( E(T)))=P_{X}C^*(E)P_{X}$.
\end{prop}

\begin{proof} For $v\in E(T)^0$ we have $P_X Q_v P_X=P_{s(\tau(v))}Q_vP_{s(\tau(v))}=Q_v$, and for
$e_u\in E(T)^1$ we have $P_X T_{e_u}P_X
=P_{s(\tau(s(e)))}T_{e_u}P_{s(\tau(u))}=T_{e_u}$. Hence
$\phi(C^*(E(T)))\subseteq P_X C^*(E)P_X$, and it remains to show
the opposite inclusion. To do this, we must show that the range of
$\phi$ contains all products $S_\mu S_\nu^*$ such that $\mu,\nu\in
E^*$, $s(\mu),s(\nu)\in X$ and $r(\mu)=r(\nu)$. Since for such
$\mu$ and $\nu$ we have $S_\mu S_\nu^*=S_\mu
S_{\tau(r(\mu))}^*S_{\tau(r(\mu))}S_\nu^*=(S_\mu
S_{\tau(r(\mu))}^*)(S_\nu S_{\tau(r(\nu))}^*)^*$, we may assume
that $\nu=\tau(r(\mu))$. The proof is by induction on the length
of $\mu$.

If $|\mu|=0$ then $\mu=s(\mu)\in X$ and so $\mu=\tau(r(\mu))$.
Then $S_\mu
S_{\tau(r(\mu))}^*=S_{\tau(r(\mu))}S_{\tau(r(\mu))}^*$, which is
in the range of $\phi$ by Lemma \ref{Q_v}. Now for $n\in\N$,
suppose $|\mu|=n$ and suppose that $S_\nu S_{\tau(r(\nu))}^*$
belongs to the range of $\phi$ for all paths $\nu$ of length
$n-1$. Let $e$ be the final edge of $\mu$, and write $\mu=\mu'e$.
Then \begin{align*}S_\mu
S_{\tau(r(\mu))}^*&=S_{\mu'}S_eS_{\tau(r(e))}^*
=S_{\mu'}P_{r(\mu')}S_eS_{\tau(r(e))}^*=S_{\mu'}(S_{\tau(r(\mu'))}^*S_{\tau(r(\mu'))})S_eS_{\tau(r(\mu))}^*
\\&=(S_{\mu'}S_{\tau(r(\mu'))}^*)(S_{\tau(r(\mu'))e}S_{\tau(r(e))}^*),\end{align*}
where $S_{\mu'}S_{\tau(r(\mu'))}^*$ belongs to the range of $\phi$
by the inductive hypothesis. If $e\in T^1$ then
$\tau(r(\mu'))e=\tau(r(e))$, and so
$S_{\tau(r(\mu'))e}S_{\tau(r(e))}^*$ belongs to the range of
$\phi$ by Lemma \ref{Q_v}. If $e$ does not belong to $T^1$, then
once again we use Lemma \ref{Q_v} to give
\begin{align*}S_{\tau(r(\mu'))e}S_{\tau(r(e))}^*&=S_{\tau(r(\mu'))e}S_{\tau(r(e))}^*(S_{\tau(r(e))}S_{\tau(r(e))}^*)\\
&=S_{\tau(s(e))e}S_{\tau(r(e))}^*\left(\sum_{u\in
E(T)^0,\,r(e)\geq_T u}Q_u\right)\\&=\sum_{u\in E(T)^0,\,r(e)\geq_T
u}T_{e_u}\end{align*} which belongs to the range of $\phi$. This
completes the proof by induction.
\end{proof}

Propositions \ref{CK}, \ref{inj} and \ref{surj} prove Theorem
\ref{main}.

\begin{example} Theorem \ref{main} is already known in the case
where $X$ is hereditary: in this case we have seen that $E(T)$ is
just the subgraph of $E$ consisting of the vertices $X$ and each
edge whose source belongs to this set. Now the isomorphism of
$C^*(E(T))$ with the corner $P_X C^*(E)P_X$ is implied by
\cite[Theorem 4.1(c)]{BPRS} and its proof.
\end{example}

\begin{example} Let $p,q,r$ be positive integers, and for each nonnegative integer $n$ let
\begin{equation*}E_n=\xygraph{{u}:@(ul,ur)"u":^{\#p}[r(2)]{v}:
@(ul,ur)"v":^{\#q}[r(2)]{w}:@(ul,ur)"w",
"u":@/_3pc/^{\#(r+np)}"w"}\end{equation*}

\noindent Example \ref{pqr} and Theorem \ref{main} show that for
each $n$ we have $C^*(E_n)\cong P_uC^*(E_0)P_u$. Graphs of this
type arise, for example, in the study of quantum lens spaces
\cite{HSQLS}.
\end{example}

\section{Labellings of directed graphs and discrete
coactions}\label{coactions}

Corners of graph algebras arise as the fixed point algebras of
certain discrete coactions on graph algebras. In this section we
make this precise, and then link up with Theorem \ref{main}.

Throughout this section, $G$ will be a group with identity element
$1_G$, or $0_G$ if $G$ is abelian (we'll sometimes omit the
subscript $G$ to avoid clutter). We will assume throughout that
$G$ is discrete, with the exception of Corollary
\ref{abelianfixedpoint} and the paragraph immediately following
the statement of Theorem \ref{fixedpoint}. We denote by $C^*(G)$
the full group $C^*$-algebra, and for $s\in G$ we simply write $s$
to denote the image of the group element under the canonical
mapping $G\to C^*(G)$. We denote by $\lambda$ the left-regular
representation of $G$ (and $C^*(G)$) on $l^2(G)$, and by $M$ the
representation of $C_0(G)$ on $l^2(G)$ by multiplication. For
$s\in G$, $\chi_s$ denotes the characteristic function of $\{s\}$.
All $C^*$-algebra tensor products considered here will involve at
least one nuclear $C^*$-algebra, and we write $A\otimes B$ to mean
the completion of the algebraic tensor product in its unique
$C^*$-norm.

\subsection*{Coactions of discrete groups on $C^*$-algebras} In line with \cite{DPR} we adopt
the following notations and conventions. Let $\delta_G:C^*(G)\to
C^*(G)\otimes C^*(G)$ be the \emph{comultiplication} $s\mapsto
s\otimes s$ for $s\in G$. A \emph{coaction} of $G$ on $A$ is an
injective nondegenerate homomorphism $\delta:A\to A\otimes C^*(G)$
(where ``nondegenerate'' means that $\lsp\delta(A)A\otimes C^*(G)$
is dense in $A\otimes C^*(G)$) such that $(\delta\otimes\id)\circ
\delta=(\id\otimes\delta_G)\circ\delta$. For each $s\in G$ let
$A_s=\{a\in A:\delta(a)=a\otimes s\}$, and write $a_s$ to denote a
generic element of $A_s$. The span of the subalgebras $A_s:s\in G$
is dense in $A$ (here it is important that $G$ be discrete). The
\emph{fixed point algebra} $A^\delta$ of $\delta$ is defined as
$A^\delta=A_{1_G}$.

A \emph{covariant representation} of the triple $(A,G,\delta)$ is
a pair $(\pi,\mu)$, where $\pi:A\to B(\mathcal H)$ and
$\mu:C_0(G)\to B(\mathcal H)$ are nondegenerate representations on
a Hilbert space $\mathcal H$ satisfying
$\pi(a_s)\mu(\chi_t)=\mu(\chi_{st})\pi(a_s)$ for all $s,t\in G$
and $a_s\in A_s$. Given a nondegenerate representation $\pi$ of
$A$ on $\mathcal H$, and letting $\lambda$ denote the left-regular
representation of $C^*(G)$ on $l^2(G)$ and $M$ the representation
of $C_0(G)$ on $l^2(G)$ by multiplication, there is a covariant
representation $((\pi\otimes\lambda)\circ\delta,1\otimes M)$ of
$(A,G,\delta)$ on $\mathcal H\otimes l^2(G)$, called the
\emph{regular covariant representation induced by} $\pi$
\cite[Proposition 2.6]{R}. The coaction $\delta$ is \emph{normal}
if there is a covariant representation $(\pi,\mu)$ with $\pi$
faithful. The \emph{crossed product} $A\times_\delta G$ is the
universal $C^*$-algebra generated by a covariant representation
$(j_A,j_G)$, and is densely spanned by elements of the form
$j_A(a_s)j_G(\chi_t)$. The nondegenerate representations of
$A\times_\delta G$ are in one-to-one correspondence with the
covariant representations of $(A,G,\delta)$ \cite[Definition
2.8]{R}. In particular, any nondegenerate representation $\pi$ of
$A$ on a Hilbert space $\mathcal H$ induces a nondegenerate
representation of the crossed product on $\mathcal H\otimes
l^2(G)$, corresponding to the regular covariant representation
induced by $\pi$. We denote this representation of $A\times_\delta
G$ by $\Ind\pi$.

\subsection*{Labellings of directed graphs} Once again, conventions are
adopted from \cite{DPR}; see also \cite{GT,KP}. A \emph{labelling}
of a directed graph $E$ by a discrete group $G$ is a function
$c:E^1\to G$. Given a labelling $c$, the \emph{skew product}
$E\times_c G$ (called the \emph{voltage graph} in \cite{GT}) is
the graph with vertex set $E^0\times G$, edge set $E^1\times G$
and whose source and range maps are given by $s(e,s)=(s(e),c(e)s)$
and $r(e,s)=(r(e),s)$ (there are several slightly different
definitions of the skew product, all yielding isomorphic graphs
\cite{DPR,GT,KQR,KP}). As further notation, if $c$ is a labelling
of a directed graph $E$ by a group $G$ and if $\mu$ is an element
of $E^*$, we shall denote by $c(\mu)$ the element
$c(\mu_1)c(\mu_2)\dots c(\mu_{|\mu|})$ of $G$.

\begin{example}\label{O_22} Let $E$ be the graph with one vertex $v$ and two
edges $e$ and $f$. Define a labelling $c$ of $E$ by $\Z_3$ as
$c(e)=2$, $c(f)=1$. Then the skew product graph $E\times_c\Z_3$ is
as follows:

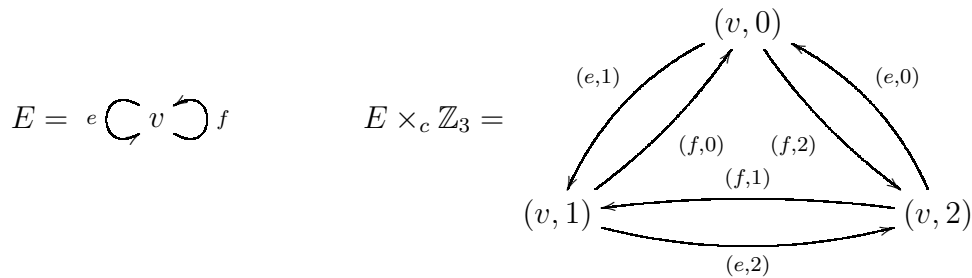
\[E=\xygraph{ v="v":@(ul,dl)_e"v":@(dr,ur)_f"v"}\qquad\qquad E\times_c\Z_3=\xygraph{[u]{(v,0)}="v":_{(e,1)}
@/_1pc/[dl(2)]{(v,1)}="dl":_{(e,2)}@/_1pc/[r(4)]{(v,2)}="dr":_{(e,0)}@/_1pc/"v":_-{(f,2)}
@/_/"dr":_-{(f,1)}@/_/"dl":_-{(f,0)}@/_/"v"   } \]
\end{example}

Any labelling of a graph $E$ by a discrete group $G$ induces a
normal coaction of $G$ on $C^*(E)$, and the graph $C^*$-algebra of
the skew product is naturally isomorphic to the crossed product of
$C^*(E)$ by this coaction; this is the content of the following
results from \cite{KQR} and \cite{DPR}, which we recall here for
convenience:

\begin{lemma}\label{KQR2.3}\cite[Lemma 2.3]{KQR}\cite[Lemma 3.3]{DPR} Let $c$ be a
labelling of a directed graph $E$ by a discrete group $G$. Then
there is a normal coaction $\delta$ of $G$ on $C^*(E)$ such that
\[\delta(S_e)=S_e\otimes c(e)\quad\text{and}\quad
\delta(P_v)=P_v\otimes 1_G\quad\text{for $e\in E^1$, $v\in
E^0$.}\]
\end{lemma}

\begin{theorem}\label{DPR3.4}\cite[Theorem 3.4]{DPR} Let $c$ be a
labelling of a directed graph $E$ by a discrete group $G$, with
corresponding coaction $\delta$. Then
\[C^*(E\times_c G)\cong C^*(E)\times_{\delta} G\] under the
isomorphism $\phi$ given on generators by
\[\phi(P_{(v,s)})=j_{C^*(E)}(P_v)j_G(\chi_s)\quad\text{and}\quad
\phi(S_{(e,s)})=j_{C^*(E)}(S_e)j_G(\chi_s)\] for $v\in E^0$, $e\in
E^1$ and $s\in G$.
\end{theorem}

\begin{example}\label{Cayley} Let $G$ be a discrete group and $S\subseteq G$ a
generating subset of cardinality $n\in\N\cup\{\infty\}$.
\cite[Theorem 2.2.3]{GT} implies that any Cayley graph $\Gamma$
for the pair $(G,S)$ is isomorphic to a skew product, by $G$, of
the graph $B_n$ with one vertex and $n$ edges. Hence by
\cite[Theorem 3.4]{DPR}, there is a normal coaction $\delta$ of
$G$ on $C^*(B_n)$ (=$\mathcal O_n$, the Cuntz algebra generated by
$n$ nonunitary isometries with mutually orthogonal ranges
\cite{C}), such that \[C^*(\Gamma)\cong \mathcal O_n\times_\delta
G.\]
\end{example}

\begin{remark} \cite[Theorem 3.4]{DPR} is more general than the
version used here. \cite{DPR} considered ``coactions of a
homogeneous space'' $G/H$, and defined an analogue of the crossed
product $A\times_\delta(G/H)$. Similarly, one may define skew
products of graphs by discrete homogeneous spaces rather than
discrete groups. \cite[Theorem 3.4]{DPR} then says that
$C^*(E\times_c(G/H))\cong C^*(E)\times_\delta(G/H)$. At this stage
it is not clear how to extend our Theorem \ref{fixedpoint} to this
more general setting, as we don't have an obvious analogue of the
fixed point algebra.\end{remark}

\subsection*{Fixed point algebras associated to labellings} We shall use the isomorphism of \cite[Theorem 3.4]{DPR}
to prove the following:

\begin{theorem}\label{fixedpoint} Let $c$ be a labelling of a directed graph $E$ by a discrete group
$G$, with corresponding coaction $\delta$. Then
\[C^*(E)^{\delta}\cong P_{E^0\times\{1_G\}}C^*(E\times_c
G)P_{E^0\times\{1_G\}}.\]\end{theorem}

This will follow immediately from the next lemma, which is the
analogue for coactions of a well-known fact about compact group
actions on $C^*$-algebras: if $\alpha:G\to\Aut(A)$ is such an
action and $i_G:G\to M(A\times_\alpha G)$ the canonical embedding,
then $i_G(1)$ is a projection which compresses the crossed product
to a copy of the fixed point algebra $A^\alpha$. Our lemma is a
weaker version of a result proved by Quigg \cite[Corollary
2.5]{Q}; we give a proof here for the sake of completeness.

\begin{lemma}\label{fixedpointlemma} Let $\delta$ be a coaction of
a discrete group $G$ on a $C^*$-algebra $A$, and let $(j_A,j_G)$
be the universal covariant representation of $(A,G,\delta)$. Then
the fixed point algebra $A^\delta$ is isomorphic to the corner
$j_G(\chi_1)(A\times_\delta G)j_G(\chi_1)$.
\end{lemma}

\begin{proof} Consider the linear map $\psi:A^\delta\to
A\times_\delta G$, $a_1\mapsto j_{A}(a_1)j_G(\chi_1)$. For
$a_1,b_1\in A^\delta$ the covariance property gives
\[j_A(a_1b_1)j_G(\chi_1)=j_A(a_1)j_A(b_1)j_G(\chi_1)^2=j_A(a_1)j_G(\chi_1)j_A(b_1)j_G(\chi_1)\]
and
\[j_A(a_1^*)j_G(\chi_1)=\left(j_G(\chi_1)j_A(a_1)\right)^*=\left(j_A(a_1)j_G(\chi_1)\right)^*,\]
so $\psi$ is a homomorphism. Furthermore, it is injective: let
$\pi$ be a faithful representation of $A$ on a Hilbert space
$\mathcal H$, and consider the induced representation $\Ind\pi$ of
$A\times_\delta G$ on $\mathcal H\otimes l^2(G)$. Then for each
nonzero $a_1\in A^\delta$ we have
\[\Ind\pi(\psi(a_1))=(\pi(a_1)\otimes\lambda_1)(1\otimes
M(\chi_1))=\pi(a_1)\otimes M(\chi_1)\neq 0.\] Hence we have shown
that $A^\delta\cong j_A(A^\delta)j_G(\chi_1)$.

For $s,t\in G$ and $a_s\in A_s$ the covariance of $(j_A,j_G)$
implies
\[
j_G(\chi_1)\left(j_A(a_s)j_G(\chi_t)\right)j_G(\chi_1)=\begin{cases}j_A(a_s)j_G(\chi_1)&\text{if
}s=t=1\\ 0&\text{otherwise.}\end{cases}\] The span of the elements
$j_A(a_s)j_G(\chi_t)$ is dense in $A\times_\delta G$, and
everything in sight is continuous, so this shows that
\[j_G(\chi_1)\left(A\times_\delta
G\right)j_G(\chi_1)=\overline{j_A(A^\delta)j_G(\chi_1)}\cong
A^\delta.\]

\end{proof}

\begin{proof}[Proof of Theorem \ref{fixedpoint}] To simplify the notation, write $P$ for the element
$P_{E^0\times \{1_G\}}$ of $ M(C^*(E\times_c G))$. The isomorphism
$\phi$ of Theorem \ref{DPR3.4} extends to a strictly continuous
isomorphism $\bar{\phi}$ of $ M(C^*(E\times_c G))$ onto $
M(C^*(E)\times_{\delta} G)$, such that
\[\bar{\phi}(P)=\sum_{v\in E^0}\phi(P_{(v,1)})=\sum_{v\in
E^0}j_{C^*(E)}(P_v)j_G(\chi_1)=j_{C^*(E)}\left(\sum_{v\in
E^0}P_v\right)j_G(\chi_{1})=j_G(\chi_{1}).\] (Note that in writing
``$j_{C^*(E)}\left(\sum_{v\in E^0}P_v\right)$'' we are actually
using the strongly continuous extension of $j_{C^*(E)}$ to the
multiplier algebra $M(C^*(E))$.) Thus the corner in which we are
interested is isomorphic to $j_G(\chi_{1})(C^*(E)\times_\delta G)
j_G(\chi_{1})$, and so Lemma \ref{fixedpointlemma} proves the
theorem.
\end{proof}

\begin{example}\label{O_23} Continuing on from Examples \ref{O_2} and \ref{O_22}, the coaction
$\delta:C^*(E)\to C^*(E)\otimes C^*(\Z_3)$ induced by the
labelling $c$ is defined on generators by $\delta(P_v)=P_v\otimes
0_{\Z_3}$, $\delta(S_e)=S_e\otimes 2$ and $\delta(S_f)=S_f\otimes
1$. Proposition \ref{fixedpoint} says that the fixed point algebra
$C^*(E)^\delta$ is isomorphic to the corner
$P_{(v,0)}C^*(E\times_c\Z_3)P_{(v,0)}$. As in Example \ref{O_2},
let $T$ be the row- and path-finite directed subtree of
$E\times_c\Z_3$ with vertices $E^0\times\Z_3$ and edges $(e,1)$
and $(f,2)$. Then Theorem \ref{main} implies that $C^*(E)^\delta$
is isomorphic to $C^*(E(T))$, where $E(T)$ is the graph
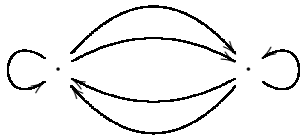
\[\xygraph{
{\cdot}="v_1":@(ul,dl)"v_1":@/^2pc/[r(2)]{\cdot}="v_2"
:@(dr,ur)"v_2":@/^2pc/"v_1":@/^1pc/"v_2":@/^1pc/"v_1"}
\]\end{example}

When the group $G$ is abelian, coactions of $G$ correspond (via
the Fourier transform) to actions of the dual group $\widehat G$
of group homomorphisms $\chi:G\to\T$ \cite[Remark 2.7]{EQ}. The
following corollary applies this fact to Theorem \ref{fixedpoint}.

\begin{cor} \label{abelianfixedpoint} Let $\alpha$ be an action of a compact abelian group $G$
on a graph algebra $C^*(E)$, such that for each $t\in G$, each
$v\in E^0$ and each $e\in E^1$, $\alpha_t(P_v)=P_v$ and
$\alpha_t(S_e)=\chi(e,t)S_e$ for some $\chi(e,t)\in\C$. Then there
is a labelling $c$ of $E$ by $\widehat G$ such that
\[C^*(E)^\alpha\cong P_{E^0\times\{1_G\}}C^*(E\times_c \widehat
G)P_{E^0\times\{1_G\}}.\]
\end{cor}

\begin{proof} Each $\alpha_t$ is an automorphism of $C^*(E)$, so each
$\chi(e,t)S_e$ is a partial isometry. This implies that for all
$e\in E^1$ and $t\in G$, $\chi(e,t)\in\T$. Now as
$\alpha:G\to\Aut(C^*(E))$ is a group homomorphism, we must have
each $\chi(e,\cdot):G\to\T$ a group homomorphism, so
$\chi(e,\cdot)\in\widehat G$ for each $e\in E^1$. Let $c$ be the
labelling $e\mapsto\chi(e,\cdot)$ of $E$ by the discrete group
$\widehat G$, and let $\delta$ be the induced coaction of
$\widehat G$ on $C^*(E)$. The Fourier transform then gives
$C^*(E)^\alpha=C^*(E)^\delta$, so the corollary follows from
Theorem \ref{fixedpoint}.
\end{proof}

The following key example was first brought to our attention by
David Pask:

\begin{example} Let $E$ be a directed graph, and let $c:E^1\to\Z$
be the labelling $c(e)=1$ for all $e\in E^1$. This labelling
corresponds to the canonical gauge action $\gamma$ of $\T$ on
$C^*(E)$, so Corollary \ref{abelianfixedpoint} implies that the
$AF$-core $C^*(E)^\gamma$ is isomorphic to the corner $
P_{E^0\times\{0_\Z\}}C^*(E\times_c\Z)P_{E^0\times\{0_\Z\}}$.
\end{example}

Theorem \ref{fixedpoint} and its corollary lead us to seek
conditions on the labelling $c$ which allow us to apply Theorem
\ref{main} to find a graph for $C^*(E)^{\delta}$. That is, we seek
conditions on $c$ which ensure that $E\times_c G$ has a row- and
path-finite directed subtree with roots $E^0\times\{1_G\}$ and
vertex set $H_{E\times_c G}(E^0\times\{1_G\})$. If $E^0$ and $G$
are both finite, the graph $E\times_c G$ has finitely many
vertices and so Lemma \ref{extree} tells us that we can always
find such a tree. More generally, when $E$ is row-finite we have
the following:

\begin{prop} \label{Kirchhoff} Let $c$ be a labelling of a row-finite directed
graph $E$ by a discrete group $G$ such that for each $\mu\in
E^\infty$, there exists $i\in\N$ such that
$c(\mu_1\ldots\mu_i)=1_G$. Then $E\times_c G$ has a row- and
path-finite directed subtree with roots $E^0\times\{1_G\}$ and
vertex set $H_{E\times_c G}(E^0\times\{1_G\})$.\end{prop}

(Note that when $E$ is finite and $G=\Z$, the above condition is
equivalent to the condition that $c(\lambda)=0$ for each loop
$\lambda\in E^*$, in analogy with Kirchhoff's voltage law.)

\begin{proof} We may apply the iterative method used in the proof
of Lemma \ref{extree} to find a directed subtree $T$ of $E\times_c
G$ with the desired roots and vertex set. This tree is row-finite
because $E$ is, so it remains to show that it is path-finite.
Suppose it is not, and let $\mu=(\mu_1,t_1)(\mu_2,t_2)\ldots$ be
an infinite path in $T$. Lemma \ref{T}(1) implies that we may
assume that the source of $\mu$ is in $E^0\times\{1_G\}$, so that
$t_1=c(\mu_1)^{-1}$, and then by definition of the source and
range maps in $E\times_c G$ we must have $r(\mu_i)=s(\mu_{i+1})$
and $t_i =c(\mu_1\ldots\mu_i)^{-1}$ for each $i\in\N$. Now
$\mu_1\mu_2\ldots$ is an infinite path in $E^\infty$, so by
assumption there exists an index $i\in\N$ such that
$c(\mu_1\ldots\mu_i)=1_G$. For this $i$ we have $r\left((\mu_i,t_i
)\right)=(r(\mu_i),c(\mu_1\ldots\mu_i)^{-1})=(r(\mu_i),1_G)$, so
the vertex $(r(\mu_i),1_G)$ receives an edge in the subtree $T$.
This is a contradiction, since each vertex in $E^0\times\{1_G\}$
is a root of $T$, and so we conclude that $T$ is path-finite.
\end{proof}

\begin{remark} We can also prove the following analogue of Proposition
\ref{Kirchhoff} for coactions of homogeneous spaces: when $E$ is
row-finite and $H$ is a subgroup of $G$ such that for each $\mu\in
E^\infty$ there exists $i\in\N$ with $c(\mu_1\ldots\mu_i)\in H$,
then there is a directed subtree in $E\times_c (G/H)$ with the
desired properties. The proof is virtually identical to the proof
of the preceding proposition, and is therefore
omitted.\end{remark}


\begin{thebibliography}{00}

\bibitem{BHRS} T. Bates, J. H. Hong, I. Raeburn and W.
Szyma\'nski, \emph{The ideal structure of the $C^*$-algebras of
infinite graphs}, Illinois J. Math. {\bf 46} (2002), 1159--1176.

\bibitem{BP} T. Bates and D. Pask, \emph{Flow equivalence of graph
algebras}, Ergodic Theory \& Dynamical Systems {\bf 24} (2004),
367--382.

\bibitem{BPRS} T. Bates, D. Pask, I. Raeburn and W.
Szyma\'nski, \emph{The $C^*$-algebras of row-finite graphs}, New
York J. Math. {\bf 6} (2000), 307--324.

\bibitem{CG} T. Crisp and D. Gow, \emph{Contractible subgraphs and
Morita equivalence of graph $C^*$-algebras}, Proc. Amer. Math.
Soc., to appear. [arXiv:math.OA/0404542]

\bibitem{C} J. Cuntz, \emph{Simple $C^*$-algebras generated by
isometries}, Commun. Math. Phys. {\bf 57} (1977), 173--185.

\bibitem{CK} J. Cuntz and W. Krieger, \emph{A class of
$C^*$-algebras and topological Markov chains}, Invent. Math. {\bf
56} (1980), 251--268.

\bibitem{DPR} K. Deicke, D. Pask and I. Raeburn, \emph{Coverings
of directed graphs and crossed products of $C^*$-algebras by
coactions of homogeneous spaces}, Internat. J. Math. {\bf 14}
(2003), 773--789.

\bibitem{D} D. Drinen, \emph{Viewing AF-algebras as graph
algebras}, Proc. Amer. Math. Soc. {\bf 128} (1999), 1991--2000.

\bibitem{DT} D. Drinen and M. Tomforde, \emph{The $C^*$-algebras of
arbitrary graphs}, Rocky Mountain J. Math. {\bf 35} (2005),
105--135.

\bibitem{DT2} D. Drinen and M. Tomforde, \emph{Computing $K$-theory
and $\operatorname{Ext}$ for graph $C^*$-algebras}, Illinois J.
Math. {\bf 46} (2002), 81--91.

\bibitem{EQ} S. Echterhoff and J. Quigg, \emph{Induced coactions
of discrete groups on $C^*$-algebras}, Canad. J. Math. {\bf 51}
(1999), 745--770.

\bibitem{FLR} N. Fowler, M. Laca and I. Raeburn, \emph{The
$C^*$-algebras of infinite graphs}, Proc. Amer. Math. Soc. {\bf
128} (2000), 2319--2327.

\bibitem{GT} J. Gross and T. Tucker, \emph{Topological graph
theory} (Wiley-Interscience, New York, 1987).

\bibitem{HSQSPSGA} J. H. Hong and W. Szyma\'nski, \emph{Quantum spheres
and projective spaces as graph algebras}, Commun. Math. Phys. {\bf
232} (2002), 157--188.

\bibitem{HSQLS} J. H. Hong and W. Szyma\'nski, \emph{Quantum lens
spaces and graph algebras}, Pacific J. Math. {\bf 211} (2003),
249--263.

\bibitem{HSPIS} J. H. Hong and W. Szyma\'nski, \emph{The primitive
ideal space of the $C^*$-algebras of infinite graphs}, J. Math.
Soc. Japan {\bf 56} (2004), 45--64.

\bibitem{KQR} S. Kaliszewski, J. Quigg and I. Raeburn, \emph{Skew
products and crossed products by coactions}, J. Operator Theory
{\bf 46} (2001), 411--433.

\bibitem{KP} A. Kumjian and D. Pask, \emph{$C^*$-algebras of
directed graphs and group actions}, Ergodic Theory \& Dynamical
Systems {\bf 19} (1999), 1503--1519.

\bibitem{KPR} A. Kumjian, D. Pask and I. Raeburn,
\emph{Cuntz-Krieger algebras of directed graphs}, Pacific J. Math.
{\bf 184} (1998), 161--174.

\bibitem{KPRR} A. Kumjian, D. Pask, I. Raeburn and J. Renault,
\emph{Graphs, groupoids and Cuntz-Krieger algebras}, J. Funct.
Anal. {\bf 144} (1997), 505--541.

\bibitem{Q} J. Quigg, \emph{Discrete $C^*$-coactions and
$C^*$-algebraic bundles}, J. Austral. Math. Soc. (Ser. A) {\bf 60}
(1996), 204--221.

\bibitem{R} I. Raeburn, \emph{On crossed products by coactions and
their representation theory}, Proc. London Math. Soc. {\bf 64}
(1992), 625--652.

\bibitem{RS} I. Raeburn and W. Szyma\'nski, \emph{Cuntz-Krieger
algebras of infinite graphs and matrices}, Trans. Amer. Math. Soc.
{\bf 358} (2004), 39--59.

\bibitem{SzK} W. Szyma\'nski, \emph{The range of $K$-invariants
for $C^*$-algebras of infinite graphs}, Indiana Univ. Math. J.
{\bf 51} (2002), 239--249.

\bibitem{SzCKU} W. Szyma\'nski, \emph{General Cuntz-Krieger uniqueness
theorem}, Internat. J. Math., {\bf 13} (2002), 549--555.

\end{thebibliography}
\end{document}